\documentclass[a4paper,11p]{article}
\usepackage{setspace}
\usepackage[english]{babel}
\usepackage[latin1]{inputenc}
\usepackage{alltt}
\usepackage[T1]{fontenc}
\usepackage{geometry}
\usepackage{graphicx,color}
\usepackage{epstopdf}
\usepackage[usenames,dvipsnames]{xcolor}
\usepackage{amsfonts,amsmath,amssymb,amsthm,lipsum}
\usepackage[normalem]{ulem}
\usepackage{eurosym}
\usepackage{hyperref}
\usepackage{listings}
\usepackage{indentfirst}
\usepackage{tikz}
\usepackage[font=small,labelfont=bf]{caption}
\usepackage{subfigure}
\usepackage{authblk}
\usepackage{pgf,tikz,pgfplots}
\pgfplotsset{compat=1.14}
\usepackage{mathrsfs}
\usepackage{pstricks-add}
\usepackage{hhline,array}
\usepackage{enumitem}
\usepackage{float}
\restylefloat{table}
\usetikzlibrary{arrows}
\definecolor{qqqqff}{rgb}{0,0,1}
\theoremstyle{break}
\newtheorem{teo}{Theorem}[section]
\newtheorem{lema}{Lemma}[section]

\newtheorem{ex}{Example}[section]
\newtheorem{prop}{Proposition}[section]
\newtheorem{defi}{Definition}[section]
\newtheorem{algo}{Algorithm}[section]
\newcommand\ChangeRT[1]{\noalign{\hrule height #1}}

\definecolor{light-gray}{gray}{0.95}
\definecolor{dkgreen}{rgb}{0,0.6,0}
\definecolor{gray}{rgb}{0.5,0.5,0.5}
\definecolor{mauve}{rgb}{0.58,0,0.82}

\title{Denoising and Interior Detection Problems}

\author[1]{Nuno Picado}
\author[1]{Paulo Eduardo Oliveira}
\affil[1]{University of Coimbra, CMUC, Department of Mathematics}
\date{}
\begin{document}
	\psset{xunit=1cm,yunit=1cm,algebraic=true,dimen=middle,dotstyle=o,dotsize=5pt 0,linewidth=1.6pt,arrowsize=3pt 2,arrowinset=0.25}
\maketitle

\begin{center}
	\textbf{Abstract}
\end{center}

Let $\mathcal{M}$ be a compact manifold of $\mathbb{R}^d$. The goal of this paper is to decide, based on a sample of points, whether the interior of $\mathcal{M}$ is empty or not. We divide this work in two main parts.
Firstly, under a dependent sample which may or may not contain some noise within, we characterize asymptotic properties of an interior detection test based on a suitable control of the dependence.
Afterwards, we drop the dependence and consider a model where the points sampled from the manifold are mixed with some points sampled from a different measure (noisy observations). We study the behaviour with respect to the amount of noisy observations, introducing a methodology to identify true manifold points, characterizing convergence properties.


\section{Introduction}

In recent years we have observed an increasing interest in estimating or testing about topological properties of an underlying set that is sampled. Indeed, in many applications, key properties of the behavior of the random phenomenon generating the observations are expressed through the geometrical or topological complexity. One is often faced with a rather high dimensional collection of observations where only a relatively small  amount of the coordinates are actually meaningful, implying that the points are, in fact, supported by, or at least are close, a lower dimensional set. The reconstruction of manifolds, based on a sample of points, has received a lot of attention, and has well established methodologies using simplicial structures. These are quite complex and computationally demanding, especially if one is interested in higher dimensional problems. However, with the help of appropriate topological tools, it is possible to obtain approximation procedures that conveniently describe the support of the observed points and some of its intrinsic complexity (such as the ones described in Medina \& Doerge \cite{medina})

The most obvious geometrical features of a set are its size, measured by the (Lebesgue) measure or some other variations on obtaining information about size, as in, for example, Pateiro-López \& Rodríguez-Casal \cite{PateiroRodriguez}, Carlstein \& Krishnamoorthy \cite{CK92}, Baíllo \& Cuevas \cite{BAC01}, Cuevas \& Rodríguez-Casal \cite{CRC04} and Cuevas, Fraiman \& Rodríguez-Casal \cite{CCF07}. Concerning the reconstruction of, possibly high dimensional, sets, we mention more elementary approaches as compared to the topologically inspired ones referred above, as these will be at the base of the results proved in this paper. A simple set estimator was introduced by Devroye and Wise~\cite{Devroye}, and, of course, variants of the same ideia have appeared elsewhere in the literature (for example, in Cuevas \cite{setestimation}). More difficult to address, from the statistical point of view, are properties such orientation or dimensionality of a set. We will contribute for a version of the later problem, considering a short procedure to decide whether a set is of full dimensionality. This will be achieved through approximating the interior of the support and deciding about is emptiness, extending the approach introduced in Aaron, Cholaquidis \& Cuevas \cite{Cuevas}, considering dependent samples, introducing and describing a suitable control so that the asymptotic characterizations introduced in \cite{Cuevas} still provide a reliable method for deciding about the interior.

As often happens in statistical problems, the observations may be subject to some noise, being interpreted here as some points in the sample possibly not being in the set of interest. As a second approach, we study a denoising procedure for selecting true points from the set, so we can afterwards apply the interior decision methods. We shall characterize how much noise is allowed in that sample so that we may still rely on the approximations for the testing procedures.

\section{Background}\label{back}

Here, we begin by introducing the reader to some geometrical definitions that will be used later, and also give some relevant results about the relation between some of them.

We first recall the definition of the estimator introduced by Devroye \& Wise \cite{Devroye}, that will at the base of most of our work.

\begin{defi}\label{estimator}
	Let $\mathcal{X}_n=\{X_1,\dots,X_n\}$ be a random sample of points in $\mathbb{R}^d$. Given $r>0$, the Devroye-Wise estimator is defined by 
	$$\hat{S}_n(r)=\bigcup_{i=1}^nB(X_i,r),$$
	where $B(x,r)$ is the closed ball with center at $x$ and radius $r$. 
\end{defi}

 In \cite{Devroye}, $\hat{S}_n(r)$ is proved to be a consistent estimator in the sense that, under suitable conditions, the measure of $\mathcal{M}\bigtriangleup\widehat{\mathcal{M}}$ will converge to $0$.

\begin{defi}
	Let $S_1,S_2\subseteq\mathbb{R}^n$. The Hausdorff distance between $S_1$ and $S_2$ is defined by
	$$d_H(S_1,S_2)=\max\{\sup_{x\in S_1}\inf_{y\in S_2}||x-y||,\sup_{y\in S_2}\inf_{x\in S_1}||x-y||\}.$$
\end{defi}

In order to be able to prove the main results of this paper, the manifold $\mathcal{M}$ will have to fulfill some regularity conditions, mainly concerning the boundary, such as the ones that are defined next.

\begin{defi}
	Let $\mathcal{M}$ be a set in $\mathbb{R}^d$. $\mathcal{M}$ is said to fulfill the outside $r$-rolling condition if for all $s\in\partial \mathcal{M}$ there exists $x\in \mathcal{M}^c$ such that $B(x,r)\cap\partial \mathcal{M}=\{s\}$. Moreover, $\mathcal{M}$ is said to fulfill the inside $r$-rolling condition if $\mathcal{M}^c$ satisfies the outside $r$-rolling condition.
\end{defi}

Intuitively, this means that by rolling a ball of radius $r$ in the border of $\mathcal{M}$, one can touch every point of this border without touching more than one at the same time. This definition obviously depends on the chosen metric.

\begin{defi}
	Let $\textrm{reach}(\mathcal{M},x)=\inf\{||x-y||:y\in\textrm{Ma}(\mathcal{M})\}$, where $\textrm{Ma}(\mathcal{M})$, called the medial axis of $\mathcal{M}$, is the set of points in $\mathbb{R}^d$ with more than one orthogonal projection onto $\mathcal{M}$.
	The reach of $\mathcal{M}$ is then defined by $\textrm{reach}(\mathcal{M})=\inf\{\textrm{reach}(\mathcal{M},x):x\in\mathcal{M}\}$.
\end{defi}

These two definitions are related by the next theorem, whose proof can be found in Cuevas, Fraiman \& Pateiro \cite{CueFraPat}.
\begin{teo}
	Let $\mathcal{M}\subset\mathbb{R}^d$ be a compact set with $\textrm{reach}(\mathcal{M})=r>0.$ Then $\mathcal{M}$ satisfies the outside $r$-rolling condition.
\end{teo}

This next definition might not be easy to understand in terms of the manifold itself. However, it expresses a regularity whose link to the properties introduced above is described next.

\begin{defi}
	Let $S$ be a set in $\mathbb{R}^d$. $S$ is said to be standard with constants $\delta$ and $\lambda$ and with respect to a Borel measure $\nu$ at a point $x\in S$ if 
	$$\nu(B(x,\varepsilon)\cap S)\geq\delta\mu_d(B(x,\varepsilon)),\ 0<\varepsilon\leq\lambda,$$
	where $\mu_d$ is the Lebesgue measure in $\mathbb{R}^d$. Moreover, we will denote by $\omega_d=\mu_d(B(x,1))$.
	
	A set is said to be standard if this holds for all $x\in S$.
\end{defi}

\begin{prop} If $\mathcal{M}$ satisfies the inside $r$-rolling condition and $P_{\mathcal{M}}$ has support $\mathcal{M}$ and a density $f$ bounded below by $f_0>0$, then $\mathcal{M}$ is standard with respect to $P_{\mathcal{M}}$ (with constants $\lambda=r$ and $\delta\leq\frac{f_0}{3}$).
\end{prop}
\bigskip

\textbf{Note:} In this paper we will use the notation $\mathring{\mathcal{M}}$ for the interior of the manifold $\mathcal{M}$ and $f_X$ (resp. $F_X$) for the density (resp. distribution) function of the random vector, or variable, $X$.

\section{Interior Detection}\label{interior}

Let $\mathcal{M}$ be a compact set. The goal of this section is to extend an interior identification procedure to handle suitably dependent samples. Naturally, we need to assume some kind of control on the dependence.

\begin{defi}\label{ai}
	Let $\mathcal{X}_n=\{X_1,\dots,X_n\}$ be a set of random variables. We call $\mathcal{X}_n$ a $\alpha_n$-almost independent sample, and denote by $\mathcal{AI}(\alpha_n)$ if
	$$\sup_{x=(x_1,\cdots,x_n)\in\mathbb{R}^d\times\cdots\times\mathbb{R}^d}\left|\frac{f_{(X_1,\dots,X_n)}(x)-f_{X_1}(x_1)\dots f_{X_n}(x_n)}{f_{X_1}(x_1)\dots f_{X_n}(x_n)}\right|\leq\alpha_n.$$
\end{defi}

Now, we check that it is possible to create models of samples $\mathcal{X}_n$ that fulfill $\mathcal{AI}(\alpha_n)$. In order to do this, we will introduce a way of creating joint distributions based on the marginals, using copula functions (see Nelsen \cite{Nelsen}).

\begin{defi}
	
	$C:[0,1]^n\longrightarrow[0,1]$ is called an $n$-copula if it is a joint distribution function with uniform marginals.
	
\end{defi}

By Sklar's Theorem (see \cite{Nelsen}) we know that for every vector $(X_1,\dots,X_n)$ there is a unique $n$-copula \textit{C} such that $$F_{(X_1,\dots,X_n)}(x_1,\dots,x_n)=C(F_{X_1}(x_1),\dots,F_{X_n}(x_n)).$$

Note: In the case of $n$-copulas, the condition on the dependence of $X_1,\dots,X_n$ for $d=1$ translates to $$\left|\frac{f_{(X_1,\dots,X_n)}-f_{X_1}\cdots f_{X_n}}{f_{X_1}\cdots f_{X_n}}\right|\leq\alpha_n\Leftrightarrow\left|\frac{\partial^n C}{\partial u_1\cdots\partial u_n}-1\right|\leq\alpha_n.$$

For $d=2$ the condition translates to
$$\sum_{k=n+1}^{2n}\sum_{(i_1,\dots,i_n)\in S_k}\left|\frac{\partial^kC}{\partial u_1^{i_1}\cdots\partial u_n^{i_n}}\right|\prod_{j:i_j=2}\left|\frac{J_F(x_j)e_1J_F(x_j)e_2}{f_0}\right|+\left|\frac{\partial^n C}{\partial u_1\cdots\partial u_n}-1\right|\leq\alpha_n,$$
where $S_k=\{(i_1,\dots,i_n)\in\{1,2\}^n:\sum_{j=1}^n i_j=k\}$ and $J_F$ is the  Jacobian matrix of the function $F$.

\begin{ex}
	
	Let $C:[0,1]^n\longrightarrow[0,1]$ be defined as $C(u_1,\dots,u_n)=\prod_{i=1}^nu_i+\prod_{i=1}^nf_i(u_i)$.
	In this case, we have
	
	$$\frac{\partial^n C}{\partial u_1\cdots\partial u_n}(u_1,\dots,u_n)=1+\prod_{i=1}^nf_i'(u_i).$$
	
	Hence,
	$$\left|\frac{\partial^n C}{\partial u_1\cdots\partial u_n}-1\right|\leq\alpha_n\Leftrightarrow\prod_{i=1}^n\left|f_i'(u_i)\right|\leq \alpha_n.$$
\end{ex}

Now that we have some dependence control, we will address the identification of the interior of a manifold considering dependent samples, with or without noise.

\subsection{Sampling without noise}\label{noiseless}

First, we will consider a noiseless model, where the sample comes from a distribution with support on the manifold $\mathcal{M}$.

In order to decide about the interior, we will use an estimator based on the one introduced in Definition \ref{estimator}.

\begin{defi}
	Let $\mathcal{X}_n=\{X_1,\dots,X_n\}$ be a random sample of points in $\mathbb{R}^d$. Given the Devroye-Wise estimator, a ball $B(X_i,r)$ will be called a boundary ball of $\hat{S}_n(r)$ if there exists $y\in\partial B(X_i,r)$ such that $y\in\partial\hat{S}_n(r)$. Then, the ``peeling'' of $\hat{S}_n(r)$ is defined as the union of all non-boundary balls, and will be denoted as $peel(\hat{S}_n(r))$.
\end{defi}

Throughout this section, we extend the results in \cite{Cuevas}, considering now dependent $\mathcal{AI}(\alpha_n)$ samples.

\begin{teo}\label{inicial}
Let $\mathcal{M}\subset\mathbb{R}^d$ be a compact non-empty set and $\mathcal{X}_n=\{X_1,\dots,X_n\}$ an $\mathcal{AI}(\alpha_n)$ sample with $\alpha_n$ such that $$\sum_{n=1}^\infty\alpha_n\epsilon_n^{-d}n^{-C(1-\epsilon_n)^d}<\infty,\text{ for some }C>1\text{ and some }\epsilon_n\rightarrow0.$$
Then,
\begin{enumerate}
	\item If $\mathring{\mathcal{M}}=\emptyset$ and $\mathcal{M}$ fulfils the outside rolling condition for some $r>0$, then $peel(\hat{S}_n(r'))=\emptyset$ for any $r'<r$.
	
	\item If $\mathring{\mathcal{M}}\neq\emptyset$, assume that there is a ball $B(x_0,r_0)\subset\mathring{\mathcal{M}}$ standard with respect to $P_X$, with constants $\delta$ and $\lambda$. Then $peel(\hat{S}_n(r_n))\neq\emptyset$ eventually a.s., with $r_n$ a sequence such that $\left(\frac{C}{\omega_d\delta}\frac{\log{n}}{n}\right)^{1/d}\leq r_n\leq\min\{r_0/2,\lambda\}.$
\end{enumerate}
\end{teo}

Proof:
\begin{enumerate}
	\item Repeat the arguments used for the proof of part (a) of Theorem 1 in \cite{Cuevas}.
	
	\item Choose $\{t_1,\dots,t_{\nu_n}\}$ such that $B(x_0,2r_n)\subseteq\bigcup_{i=1}^{\nu_n}B(t_i,r_n\epsilon_n)$. It is easily seen that $\nu_n=\tau_d\epsilon_n^{-d}$. Let us define $p_n:=P(\exists y\in B(x_0,2r_n),\mathring{B}(y,r_n)\cap\mathcal{X}_n=\emptyset).$
	
	Then,
	\begin{eqnarray}
		\nonumber p_n & \leq & \sum_{i=1}^{\nu_n}P(X_1\notin \underbrace{B(t_i,r_n(1-\epsilon_n))}_{B_{n,i}},\ldots,X_n\notin \underbrace{B(t_i,r_n(1-\epsilon_n))}_{B_{n,i}})\\
		\nonumber & = & \sum_{i=1}^{\nu_n}\int_{B_{n,i}^c}\cdots\int_{B_{n,i}^c} f_{(X_1,\ldots,X_n)}dx_1\ldots dx_n \\
		\nonumber& \leq & \sum_{i=1}^{\nu_n}\int_{B_{n,i}^c}\cdots\int_{B_{n,i}^c} (1+\alpha_n)f_{X_1}\cdots f_{X_n} dx_1\ldots dx_n \\
		& = & \sum_{i=1}^{\nu_n}(1+\alpha_n)P(X_1\notin B(t_i,r_n(1-\epsilon_n)))^n.\label{pn}
	\end{eqnarray}
	Note that because $r_n\leq r_0/2$, it follows that $t_i\in B(x_0,r_0)$ for $i=1,\dots,\tau_d$, so we can use the standardness of the ball $B(x_0,r_0)$ to get that
	\begin{eqnarray}
	\nonumber P(X_1\notin B(t_i,r_n(1-\epsilon_n)))^n&=&(1-P_X( B(t_i,r_n(1-\epsilon_n))))^n\\
	\nonumber&\leq&(1-\omega_d\delta r_n^d(1-\epsilon_n)^d)^n\\
	\nonumber&\leq&\left(1-C\frac{\log n}{n} (1-\epsilon_n)^d\right)^n\\
	\nonumber&\leq&\exp(-C(1-\epsilon_n)^d)\log n\\
	\nonumber&=&n^{-C(1-\epsilon_n)^n}.
	\end{eqnarray}
	
	Applying this upper bound in (\ref{pn}), we finally conclude that $$p_n\leq\tau_d\epsilon_n^{-d}(1+\alpha_n)n^{-C(1-\epsilon_n)^n}\ \Rightarrow\sum_{n=1}^{+\infty} p_n<\infty.$$
	
	By the Borel-Cantelli Lemma, it follows that for all $y\in B(x_0,2r_n)$, we have $\mathring{B}(y,r_n)\cap\{X_1,\dots,X_n\}\neq\emptyset$, eventually a.s. Putting $y=x_0$ means that there exists some $X_i\in\mathring{B}(x_0,r_n)$. Moreover, if $z\in\partial B(X_i,r_n)$ there exists $X_j$ such that $z\in\mathring B(X_j,r_n)$. This implies that the ball $B(X_i,r_n)$ belongs to $peel(\hat{S}_n(r_n))$, and so $peel(\hat{S}_n(r_n))\neq\emptyset$.
\end{enumerate}
	\begin{flushright}
		$\Box$
	\end{flushright}

The next step is to choose the radii of the balls appropriately. According to the previous theorem, these must be chosen converging to zero, but at a suitable rate. The next lemma will help us with this tuning of the decrease rate.

\begin{lema}\label{rn>}
	For $n$ large enough and $t<\frac{1}{f_0}$ for some fixed $f_0>0$, let $\{X_1,\dots,X_n\}$ be an $\mathcal{AI}(\alpha_n)$ random sample with $\alpha_n$ such that $$\sum_{n=1}^\infty\alpha_n\exp\left(-\sigma_n\frac{\gamma}{2}n^{-\tau}\right)<+\infty,$$ for any 
	$\frac{1}{2}<\tau<1$, $0<\gamma<\left(1-\left(\frac{1}{2}\right)^{1/d}\right)^d$ and $\sigma_n\geq\frac{cn}{\log n}$. 
	
	Then $\max_{i}\min_{j\neq i}||X_j-X_i||>\left(\frac{t\log{n}}{n\omega_d}\right)^{1/d}$ with probability one.
\end{lema}

Proof: Let $g_0>f_0$, $v>u>t$, $\varepsilon>0$, be such that $g_0v<1$ and $\varepsilon^{1/d}+t^{1/d}<u^{1/d}$. Let $B_0$ be a ball centered in a point $x_0\in\mathcal{M}$ and radius $\beta_0$ such that for every $x\in B(x_0,2\beta_0)$, $f(x)<g_0$. Finally let $N(n)$ and $M(n)$ be independent Poisson variables with means $n-n^{3/4}$ and $2n^{3/4}$ resp., and $\mathcal{W}_n^-=\{X_1,\dots,X_{N(n)}\}$, $\mathcal{W}_n^+=\{X_1,\dots,X_{N(n)+M(n)}\}$. Then, $\mathcal{W}_n^-$, resp. $\mathcal{W}_n^+$, is a Poisson process with intensity funcion $(n-n^{3/4})f(\cdot)$, resp. $(n+n^{3/4})f(\cdot)$ (see Kingman \cite{Kingman} for more details on Poisson processes).

With $H_n=\{\mathcal{W}_n^-\subset\mathcal{X}_n\subset\mathcal{W}_n^+\}=\{N(n)<n<N(n)+M(n)\}$, we get, by Lemma 3.1 in Penrose \cite{Penrose}, that there exists $c_1>0$ satisfying $$P(H_n^c)\leq2c_1e^{-n^{1/4}}.$$

Let $\sigma_n=\sigma(B_0,\rho_n(u))$, where $\sigma(U,r)=\max\{n:\exists x_1,\dots,x_n\in U:B(x_i,r)\cap B(x_j,r)=\emptyset\}$. By Lemma 2.1 in \cite{Penrose}, we have that $\sigma_n\geq\frac{c_2n}{\log{n}}$, for some $c_2>0$.
\medskip

Let $\{x_1^n,\dots,x_{\sigma_n}^n\}\subset B_0$ such that $B(x_i^n,\rho_n(u))\cap B(x_j^n,\rho_n(u))=\emptyset$. Given a point process $\mathcal{W}$, we denote by $\mathcal{W}[U]$ the number of points of $\mathcal{W}$ in $U$, and let $E_n(x)=\{\mathcal{P}_n^-[B(x,\rho_n(\varepsilon))]=1\}\cap\{\mathcal{P}_n^+[B(x,\rho_n(u))]=1\}$. Using that $\left\{\max_i\min_{j\neq i}||X_j-X_i||\leq\rho_n(t)\right\}\subset H_n^c\cup\left(\bigcup_{i=1}^{\sigma_n}E_n(x_i^n)\right)^c$, we get that
	
	$\begin{aligned}
	\sum_{n=1}^\infty P\left(\left\{\max_i\min_{j\neq i}||X_j-X_i||\leq\rho_n(t)\right\}\right)&\leq \sum_{n=1}^\infty P(H_n^c)+\sum_{n=1}^\infty P\left(\left(\bigcup_{i=1}^{\sigma_n}E_n(x_i^n)\right)^c\right)\\
	&\leq\sum_{n=1}^\infty2c_1e^{-n^{1/4}}+\sum_{n=1}^\infty(1+\alpha_n)\exp\left(-\sigma_n\frac{f_0\varepsilon}{2}n^{-g_0v}\right)<+\infty.
	\end{aligned}$
	
	Applying the Borel-Cantelli Lemma the result follows.

	\begin{flushright}
	$\Box$
\end{flushright}

\begin{teo}\label{teonoiseless}
	Let $\mathcal{M}$ be a $d'$-dimensional compact manifold in $\mathbb{R}^d$ and $\mathcal{X}_n$ an $\mathcal{AI}(\alpha_n)$ sample with $\alpha_n$ in the conditions of Theorem \ref{inicial} and Lemma \ref{rn>}. Assume $f(x)>f_0$ for every $x\in\mathcal{M}$. Let $r_n=\beta\max_i\min_{j\neq i}||X_j-X_i||$, with $\beta>6^{1/d}$.
	Then,
\begin{enumerate}
	\item If $d'=d$ and $\partial\mathcal{M}$ is a $\mathcal{C}^2$ manifold, then $peel(\hat{S}_n(r_n))\neq\emptyset$ eventually, a.s.
	
	\item if $d'<d$ and $\mathcal{M}$ is a $\mathcal{C}^2$ manifold without boundary, then $peel(\hat{S}_n(r_n))=\emptyset$ eventually, a.s.
\end{enumerate}
\end{teo}

Proof:
\begin{enumerate}
	\item Given the conditions ($d'=d$), we know that $\partial\mathcal{M}$ is a $\mathcal{C}^2$ compact manifold of dimension $d-1$ (see Conlon \cite{Conlon}). Therefore, by Theorem 1 in Walther \cite{Walther}, $\mathcal{M}$ fulfils the inside and outside rolling ball conditions for some $r>0$. So, by Proposition 1 in \cite{Cuevas}, $\mathcal{M}$ satisfies the standardness condition. Using Theorem \ref{inicial}, it remains to be proved that $r_n\geq \left(\kappa\frac{\log{n}}{n}\right)^{1/d}$,  for $n$ large enough and $\kappa>(\delta\omega_d)^{-1}$.
	Using Lemma \ref{rn>} with $t=\frac{1}{2f_0}$, we get that $r_n\geq\left(\beta^d\frac{\log n}{\omega_d2f_0n}\right)^{1/d}$.
	Morevover $\kappa:=\frac{\beta^d}{2\omega_df_0}>\frac{3}{2\omega_df_0}=(\omega_d\delta)^{-1}.$

	\item Due to the fact that $\mathcal{M}$ is a $\mathcal{C}^2$ compact manifold of $\mathbb{R}^d$, by Proposition 14 in Thäle \cite{Thale} it has a positive reach, implying that it satisfies the outside rolling ball condition for some $r>0$. Therefore, we may apply Theorem \ref{inicial}.
	Now, we just need to prove that $r_n\leq r$ for $n$ large enough. For that, it is enough to prove that $\max_i\min_{j\neq i}||X_j-X_i||\xrightarrow{a.s.}0$, which may be achieved reproducing the arguments as in the proof of Theorem 5.1 in \cite{Penrose}.
		
\end{enumerate}
	\begin{flushright}
		$\Box$
	\end{flushright}
	
\subsection{Sampling with noise}\label{noisy}
	
In this subsection, we will be considering the case where the sample is observed with some general random noise. We will study the case where the sample is generated from a distribution with support $\mathcal{S}=B(\mathcal{M},R)$ with a density function bounded below by $f_0>0$.
We will extend the corresponding results proved in Aaron, Cholaquidis \& Cuevas \cite{Cuevas}. 
	
\begin{teo}\label{noise}
	Let $\mathcal{M}$ be a compact set in $\mathbb{R}^d$ such that $reach(\mathcal{M})=R_0>0$. Let $\mathcal{X}_n$ be an $\mathcal{AI}(\alpha_n)$ sample of a distribution with support $\mathcal{S}=B(\mathcal{M},R_1)$ with $0<R_1<R_0$, with density $f$ bounded below by $f_0>0$ and $\alpha_n$ such that
	$$\sum_{n=1}^{+\infty}\frac{\alpha_n}{n^\beta\log n}<+\infty,\ for\ any\ \beta>1.$$ Let $\rho_n=c\left(\frac{\log{n}}{n}\right)^{1/d}$, with $c>\left(\frac{6}{f_0\omega_d}\right)^{1/d}$, $\hat{R}_n=\max_{i}\min_{j\in I_{bb}}||X_i-X_j||$ where $I_{bb}=\{j:B(Y_j,\rho_n) \mathrm{\;is\;a\;boundary\;ball}\}$.
	
	\begin{enumerate}
		\item If $\mathring{\mathcal{M}}=\emptyset$, then, with probability one, $$\left|\hat{R}_n-R_1\right|\leq 2\rho_n\;for\;n\;large\;enough,$$
		
		\item If $\mathring{\mathcal{M}}\neq\emptyset$, then there exists $C>0$ such that, with probability one $$\left|\hat{R}_n-R_1\right|>C\;for\;n\;large\;enough.$$
	\end{enumerate}
\end{teo}
	
	The proof will be presented later, after some auxiliary results.
	
\begin{lema}\label{dH}
	Let $\mathcal{X}_n$ be an $\mathcal{AI}(\alpha_n)$ sample with $\alpha_n$ such that
	$$\sum_{n=1}^{+\infty}\frac{\alpha_n}{n^\beta\log n}<+\infty,\ for\ any\ \beta>1.$$ Then,
	$$\limsup\left(\frac{n}{\log n}\right)^{1/d}d_H(\mathcal{X}_n,S)\leq\left(\frac{2}{\delta\omega_d}\right)^{1/d}\;a.s.$$	
\end{lema}

Proof: Given that $\mathcal{X}_n\subset S$, by the definition of the Hausdorff distance we know that $d_H(\mathcal{X}_n,S)=\sup_{x\in S}\min_i||x-X_i||$. Covering $S$ with balls of radius $\Delta$ and denoting by $S_\Delta$ the set of the balls centers, for every $x\in S$ and $s_0\in S_\Delta$ we find an upper bound for $\min||x-X_i||$:
\begin{align}
\nonumber\min||x-X_i||&\leq||x-X_j||\leq||x-s_0||+||s_0-X_j||\\
\nonumber&\leq||x-s_0||+\max_{s}\min_j||s-X_j||\leq\Delta+\max_{s}\min_j||s-X_j||.
\end{align}

Using $\Delta=(1-\upsilon)\varepsilon$, we get
$P(d_H(\mathcal{X}_n,S)>\varepsilon)\leq P(\max\min||s-X_i||>\upsilon\varepsilon)$. Setting $I_{s,i}=\{X_i\in B(s,\upsilon\varepsilon)\}$, we get
\begin{align}
\nonumber	P(\max\min||s-X_i||>\upsilon\varepsilon)&=P(\cup_{s}\cap_{i=1}^nI_{s,i}^c)\leq\sum_{s}P(\cap_{i=1}^nI_{s,i}^c)\\
\nonumber	&\leq(1+\alpha_n)\sum_{s}\prod_{i=1}^nP(I_{s,i}^c)\\
\nonumber	&\leq(1+\alpha_n)\sum_{s}\prod_{i=1}^n(1-P(I_{s,i})).
\end{align}

Using the standardness of the set, $P(I_{s,i})=P_X(B(s,\upsilon\varepsilon))\geq\delta\omega_d(\upsilon\varepsilon)^d$.
Then,
\begin{align}
\nonumber P(\max\min||s-X_i||>\upsilon\varepsilon)&\leq(1+\alpha_n)\sum_{s}(1-\delta\omega_d(\upsilon\varepsilon)^d)^n\\
\nonumber&\leq A(1+\alpha_n)((1-\upsilon)\varepsilon)^{-d}e^{-n\delta\omega_d(\upsilon\varepsilon)^d},
\end{align}
with $A$ a constant not depending on $n$ or $d$.

Hence, $P(d_H(\mathcal{X}_n,S)>\varepsilon)\leq A(1+\alpha_n)((1-\upsilon)\varepsilon)^{-d}e^{-n\delta\omega_d(\upsilon\varepsilon)^d}$ and

$P\left(\left(\frac{n}{\log n}\right)^{1/d}d_H(\mathcal{X}_n,S)>\ell\right)\leq A((1-\upsilon)\ell)^{-d}(1+\alpha_n)\frac{n^{1-\delta\omega_d\upsilon^d\ell^d}}{\log n}$.

The result now follows from the Borel-Cantelli lemma, using $\ell=\frac{2}{\delta\omega_d}$.

\begin{flushright}
	$\Box$
\end{flushright}

Proof of Theorem \ref{noise}: 

Since $S$ has a Lebesgue null boundary, $\hat{S}_n(\rho_n)\subset B(\mathring{S},\rho_n)$.
Due to the fact that $c>\left(\frac{6}{f_0\omega_d}\right)^{1/d}=\left(\frac{6}{3\delta\omega_d}\right)^{1/d}=\left(\frac{2}{\delta\omega_d}\right)^{1/d}$, using Lemma \ref{dH}, we get that

$$d_H(\mathcal{X}_n,S)\leq c\left(\frac{\log n}{n}\right)^{1/d}=\varepsilon_n.$$

Therefore, we conclude that, with probability one, $$S\subset\hat{S}_n(\rho_n).$$

The rest of the proof follows the same arguments as in the independent sample case.

\begin{flushright}
	$\Box$
\end{flushright}

This last approach has one really big problem when it comes to using in real data, being that we need to know beforehand the value of $R_1$. This means that we need to have some information about the amount of noise in the sample. In order to overcome this, we will introduce a different type of noise and introduce a new methodology to denoise the sample so that we can use the methodology described in subsection \ref{noiseless}.

\section{Denoising}

In this section, we will consider an independent sample of points from a probability measure $\mu_n'$ which is a mixture of two probability measures: $\mu$, whose support is the manifold $\mathcal{M}$ (with proportion $1-\alpha_n$), and $\mu_R$, considering this to be a uniform measure in a ball of radius R containing $\mathcal{M}$ (with proportion $\alpha_n$), that is, $\mu_n'=(1-\alpha_n)\mu+\alpha_n\mu_R$. The goal is to construct a method to eliminate the points that come from the second measure, therefore keeping only the points belonging to $\mathcal{M}$.

In order to do that, we will need some notions of distance between measures and distance to a measure. For the first we will use the classical Wasserstein distance (see Villani, \cite{villani}, for a more complete background).

\begin{defi}
	Given two probability measures $\mu$ and $\nu$ in $\mathbb{R}^d$, a transport plan is a probability measure $\pi$ in $\mathbb{R}^d\times\mathbb{R}^d$ s.t. $\pi(A\times\mathbb{R}^d)=\mu(A)$ and $\pi(\mathbb{R}^d\times B)=\nu(B)$.
	
	The cost of $\pi$ is defined as
	$$\mathcal{C}(\pi)=\left(\int_{\mathbb{R}^n\times\mathbb{R}^n}\|x-y\|^2d\pi(x,y)\right)^\frac{1}{2}.$$
	
	Moreover, the Wasserstein distance between two probability measures $\mu$ and $\nu$, denoted by $W_2(\mu,\nu)$ is given by the minimum of all the transport plans costs.
\end{defi}

This distance function is a really good way to quantify the distance between measures with the type of noise we have. In fact, if we consider a measure $\mu$ uniform on the set $supp(\mu)=\{x_1,\cdots,x_n\}$ and a measure $\nu$ uniform on $supp(\nu)=\{y_1,\cdots,y_k,x_{k+1},\cdots,x_n\}$, with points $y_j$ such that $\min_i||x_i-y_j||\leq R,$ for every $j\in\{1,\cdots,k\}$, we have
\begin{eqnarray}
\notag W_2(\mu,\nu)&\leq&\left(\sum_{i=1}^k\frac{1}{n}\|x_i-y_i\|^2\right)^{\frac{1}{2}}\\
\notag&\leq&\left(\sum_{i=1}^k\frac{1}{n}(R+\emph{diam}(supp(\mu)))^2\right)^{\frac{1}{2}}\\
\notag&=&\left(\frac{k}{n}(R+\emph{diam}(supp(\mu)))^2\right)^{\frac{1}{2}}\\
\notag&=&\left(\frac{k}{n}\right)^{\frac{1}{2}}(R+\emph{diam}(supp(\mu)))
\end{eqnarray}

If we consider $k<<n$, this distance will be close to zero, as we would want to happen given we are just inserting some noise in a small number of points of the sample. As for the notion of distance to a certain measure we will follow use a function introduced in Chazal, Cohen-Steiner \& Mérigot \cite{ChazalSteiner}.
\begin{defi}
	Let $\mu$ be a probability measure and $0\leq m<1$. We denote by $\delta_{\mu,m}$ the function
	$$\begin{array}{cccl}
	\delta_{\mu,m}\;: & \mathbb{R}^n & \rightarrow & \mathbb{R}_0^+ \\
	& x & \rightsquigarrow & \inf\{r>0:\mu(\bar{B}(x,r))>m\} \\
	\end{array}$$
\end{defi}

Note that in the case $m=0$, this distance coincides with the distance to the support of $\mu$.
However, this function is not robust with regard to small perturbations on the measure $\mu$. For example, define $\mu_\varepsilon=(\frac{1}{2}-\varepsilon)\delta_0+(\frac{1}{2}+\varepsilon)\delta_1$. In this case, for $\varepsilon>0$, we get $\delta_{\mu_\varepsilon,1/2}(t)=|1-t|$, while for $\varepsilon<0$ we get $\delta_{\mu_\varepsilon,1/2}(t)=|t|$. This is a problem for our methodology, as the measure we consider will change with the sample size and we need to control the distance function.

To overcome this problem, we shall consider a smoothed version of $\delta_{\mu,m}$:

\begin{defi}\label{10.5}
	Let $\mu$ be a probability measure in $\mathbb{R}^d$ and $0<m_0\leq1$. The \textbf{distance function} to $\mu$ is given by the function
	$$\begin{array}{cccl}
	d_{\mu,m_0}\;: & \mathbb{R}^n & \rightarrow & \mathbb{R}_0^+ \\
	& x & \rightsquigarrow & \left(\frac{1}{m_0}\int_0^{m_0}\delta_{\mu,m}(x)^2dm\right)^{1/2} \\
	\end{array}$$
\end{defi}

In general, this function is difficult to compute, but in the case where $\mu$ is an empirical measure it becomes much easier, as shown by this next example.

\begin{ex}\label{dm}
	Let $P$ be a set consisting of $n$ points and $\mu_P=\frac{1}{n}\sum_{p\in P} \delta_p$.
	
	Assuming $m_0=\frac{k_0}{n}$ we have
	\begin{gather}\label{d}
	d_{\mu,m_0}(x)=\left(\frac{1}{k_0}\sum_{p\in NN_P^{k_0}(x)}\|p-x\|^2\right)^{1/2},
	\end{gather}
	
	where $NN_P^{k_0}(x)$ is the set of the $k_0$ nearest neighbours of $x$ in $P$.
\end{ex}

With this new distance function, we achieve the robustness that $\delta_{\mu,m}$ lacked, as we can see by the following theorem whose proof can be found in Boissonnat, Chazal \& Yvinec \cite{Chazal}.

\begin{teo}
	Let $\mu$ and $\mu'$ be probability measures. Then $$\|d_{\mu,m_0}-d_{\mu',m_0}\|_\infty\leq m_0^{-1/2}W_2(\mu,\mu').$$
\end{teo}

\begin{teo}\label{comblin}
	Let $\mu,\mu_1,\mu_2$ be measures and define $\mu'=(1-\alpha)\mu_1+\alpha\mu_2$. Then, $$W_2(\mu,\mu')^2\leq (1-\alpha)W_2(\mu,\mu_1)^2+\alpha W_2(\mu,\mu_2)^2.$$
\end{teo}

Proof:
Let $\pi_1$ be a transport plan from $\mu$ to $\mu_1$ and $\pi_2$ a transport plan from $\mu$ to $\mu_2$.

Defining $\pi'=(1-\alpha)\pi_1+\alpha\pi_2$ we obtain a transport plan from $\mu$ to $\mu'$:

$\pi'(A\times\mathbb{R}^d)=(1-\alpha)\pi_1(A\times\mathbb{R}^d)+\alpha\pi_2(A\times\mathbb{R}^d)=(1-\alpha)\mu(A)+\alpha\mu(A)=\mu(A),$

$\pi'(\mathbb{R}^d\times B)=(1-\alpha)\pi_1(\mathbb{R}^d\times B)+\alpha\pi_2(\mathbb{R}\times B)=(1-\alpha)\mu_1(B)+\alpha\mu_2(B)=\mu'(B).$

Now,
\begin{eqnarray}
\nonumber W_2(\mu,\mu')^2 & = & \min_{\pi'}\left(\int_{\mathbb{R}^d\times\mathbb{R}^d}||x-y||^2d\pi'(x,y)\right)\\
\nonumber & \leq & \min_{\pi_1,\pi_2}\left(\int_{\mathbb{R}^d\times\mathbb{R}^d}||x-y||^2d((1-\alpha)\pi_1+\alpha\pi_2)(x,y)\right)\\
\nonumber & = & \min_{\pi_1}\left(\int_{\mathbb{R}^d\times\mathbb{R}^d}||x-y||^2d(1-\alpha)\pi_1(x,y)\right)+\min_{\pi_2}\left(\int_{\mathbb{R}^d\times\mathbb{R}^d}||x-y||^2d\alpha\pi_2(x,y)\right)\\
\nonumber & = & (1-\alpha)W_2(\mu,\mu_1)^2+\alpha W_2(\mu,\mu_2)^2.
\end{eqnarray}

\begin{flushright}
	$\Box$
\end{flushright}

\begin{teo}\label{up} Let $\mu_n'=(1-\alpha_n)\mu+\alpha_n\mu_R$ and $\widehat{\mu_n}'$ the empirical measure associated with a sample drawn from this measure. Then
	$$||d_{\mathcal{M}}-d_{\widehat{\mu_n}',m_n}||_{\infty}\leq C(\mu)^{1/d'}m_n^{1/d'}+\left(\frac{\alpha_n}{m_n}\right)^{1/2}W_2(\mu,\mu_R)+m_n^{-1/2}W_2(\mu_n',\widehat{\mu_n}'),$$
	where $C(\mu)$ is a constant depending only on $\mu$.
\end{teo}

Proof:
Using Theorem 3.5 and Corollary 4.8 in Chazal, Cohen-Steiner \& Mérigot \cite{ChazalSteiner}, and Theorem \ref{comblin} we obtain the following inequalities:

\begin{align}
\nonumber ||d_{\mathcal{M}}-d_{\widehat{\mu_n}',m_n}||_{\infty} & \leq  C(\mu)^{1/d'}m_n^{1/d'}+m_n^{-1/2}W_2(\mu,\widehat{\mu_n}')\\
\nonumber & \leq  C(\mu)^{1/d'}m_n^{1/d'}+m_n^{-1/2}W_2(\mu,\mu_n')+m_n^{-1/2}W_2(\mu_n',\widehat{\mu_n}') \\
\nonumber& \leq  C(\mu)^{1/d'}m_n^{1/d'}+m_n^{-1/2}\alpha_n^{1/2}W_2(\mu,\mu_R)+m_n^{-1/2}W_2(\mu_n',\widehat{\mu_n}')\\
& =  C(\mu)^{1/d'}m_n^{1/d'}+\left(\frac{\alpha_n}{m_n}\right)^{1/2}W_2(\mu,\mu_R)+m_n^{-1/2}W_2(\mu_n',\hat{\mu_n'})\\
& \leq  C(\mu)^{1/d'}m_n^{1/d'}+\left(\frac{\alpha_n}{m_n}\right)^{1/2}W_2(\mu,\mu_R)+m_n^{-1/2}W_2(\mu_n',\widehat{\mu_n}').
\end{align}

\begin{flushright}
	$\Box$
\end{flushright}

\begin{teo}\label{dconv} Under the same conditions as in Theorem \ref{up}, if $d\geq4$, $m_n\longrightarrow0$, $\alpha_nm_n^{-1}\longrightarrow0$, $n^{-1/d_M^*(\mu_R)}m_n^{-1/2}\longrightarrow0$, then	$d_{\hat{\mu}',m_n}\xrightarrow[]{\text{P}} d_{\mathcal{M}}$.
\end{teo}

Proof: Using Theorem \ref{up}, the only thing left to prove is that $m_n^{-1/2}W_2(\mu_n',\hat{\mu_n'})\xrightarrow[]{\text{P}} 0$. For that we use Theorem 1 in Bach \& Weed \cite{bachweed}, which states that for any $s>d_p^*(\mu_n')$, $$E\left[W_2(\mu_n',\hat{\mu_n'})\right]\leq n^{-1/s}.$$

The only problem now is finding an upper bound for $d_p^*(\mu_n')$, which is given by the fact that $d_p^*(\mu_n')\leq d_M(\mu_n')=d$ and can be proved by repeating the arguments in the proof of Proposition 2 in \cite{bachweed}

\begin{flushright}
	$\Box$
\end{flushright}

Based on Theorem \ref{dconv}, we will now introduce the method to de-noise the sample, which is basically remove all the points where $d_{\hat{\mu}',m_n}(X_i)>\delta_n$. Using theorems \ref{up} and \ref{dconv} and their proofs, we know that we will remove points where $d_{\mathcal{M}}>r_n$, with $r_n=\delta_n+C(\mu)^{1/d'}m_n^{1/d'}+\left(\frac{\alpha_n}{m_n}\right)^{1/2}W_2(\mu,\mu_R)+m_n^{-1/2}n^{-1/d}$.

\begin{teo}
	Let $\mathcal{M}\subset\mathbb{R}^d$ be a $d'$-dimensional manifold and let $m_n\sim n^{-x}$, $\alpha_n\sim n^{-y}$, $\delta_n\sim n^{-z}$.
	
	Under the conditions $$\left\{\begin{array}{l}
	1-y-x\left(\frac{d-d'}{d'}\right)<0\\
	1-y+\frac{x-y}{2}(d-d')<0\\
	1-y+\left(\frac{x}{2}-\frac{1}{d}\right)(d-d')<0\\
	1-y-z(d-d')<0
	\end{array}\right.,$$ the probability of eliminating the points not belonging to $\mathcal{M}$ will converge to 1.
\end{teo}

Proof: As we are eliminating all the points where $d_{\mathcal{M}}(X_i)>r_n$, we just need to worry about points in $B(\mathcal{M},r_n)\backslash\mathcal{M}.$
\begin{eqnarray}
\nonumber P(X_i\notin B(\mathcal{M},r_n)\backslash\mathcal{M}) & = & 1-P(X_i\in B(\mathcal{M},r_n)\backslash\mathcal{M})\\
\nonumber& = & 1-\mu_n'( B(\mathcal{M},r_n)\backslash\mathcal{M})\\
\nonumber& = & 1-(1-\alpha_n)\mu( B(\mathcal{M},r_n)\backslash\mathcal{M})-\alpha_n\mu_R( B(\mathcal{M},r_n)\backslash\mathcal{M})\\
\nonumber& = & 1-\alpha_n\mu_R( B(\mathcal{M},r_n)/\mathcal{M})\\
\nonumber& \geq & 1-c\alpha_nr_n^{d-d'}.
\end{eqnarray}

Using the independence of the sample, we need to prove that 

$$(1-c\alpha_nr_n^{d-d'} )^{n}\longrightarrow1\Leftrightarrow\frac{\log\left(1-c\alpha_nr_n^{d-d'}\right)}{n^{-1}}\longrightarrow0$$

Applying L'Hôpital's rule, this is the same as proving that

\begin{align}
\nonumber & \frac{(d-d')c\alpha_nr_n'r_n^{d-d'-1}+c\alpha_n'r_n^{d-d'}}{n^{-2}(1-c\alpha_nr_n^{d-d'})} \longrightarrow 0\\
\nonumber \Leftrightarrow\;\; & (d-d')cn^2\alpha_nr_n'r_n^{d-d'-1}+cn^2\alpha_n'r_n^{d-d'} \longrightarrow 0\\
\nonumber \Leftrightarrow\;\; & (d-d')cn^{2-y}r_n'r_n^{d-d'-1}+cn^{1-y}r_n^{d-d'} \longrightarrow 0\\
\nonumber \Leftrightarrow\;\; & n^{1-y}r_n^{d-d'-1}(nr_n'+r_n) \longrightarrow 0\\
\nonumber \Leftrightarrow\;\; & n^{1-y}r_n^{d-d'} \longrightarrow 0 .
\end{align}

Now, because $n^{1-y}r_n^{d-d'}\sim n^{1-y-x(\frac{d-d'}{d'})}+n^{1-y+\frac{x-y}{2}(d-d') }+n^{1-y+(\frac{x}{2}-\frac{1}{d})(d-d')}+n^{1-y-z(d-d')},$ using the conditions imposed in the theorem we get the result.

\begin{flushright}
	$\Box$
\end{flushright}

\begin{ex}
	If for example, we have $d=4$ and $d'=1$, we have to impose the conditions
	
	\begin{table}[H]
		\begin{center}
			\begin{tabular}{lcr}
				$\left\{\begin{array}{l}
				0<x<\frac{1}{2}\\
				y>x\\
				1-y-3x<0\\
				1-\frac{5}{2}y+\frac{3}{2}x<0\\
				\frac{1}{4}-y+\frac{3}{2}x<0
				\end{array}\right.$
				&
				$\Leftrightarrow$
				&
				\raisebox{-.5\height}{
					\begin{tikzpicture}
					\foreach \x in {0,0.25,0.5}
					\draw (4*\x cm,1pt) -- (4*\x cm,-1pt) node[anchor=north] {$\x$};
					\foreach \y in {0,0.5,1,1.5,2}
					\draw (1pt,\y cm) -- (-1pt,\y cm) node[anchor=east] {$\y$};
					\fill[blue!40!white] (0,1) -- (0.666,0.5) -- (2,1) -- (2,2.8) -- (0,2.8) -- (0,1);
					\draw[line width=0.25mm, dashed, blue!70!white] (0,2.9) -- (0,1) -- (0.666,0.5) -- (2,1) -- (2,2.9);
					\draw[thick,->] (0,0) -- (2.5,0) node[anchor=north west] {$x$};
					\draw[thick,->] (0,0) -- (0,2.5) node[anchor=south east] {$y$};
					\end{tikzpicture}}
			\end{tabular}
		\end{center}
	\end{table}
	
	along with the condition $1-y-3z<0$.
\end{ex}

We may now describe an algorithm to denoise the sample and consequently decide whether the interior of the manifold is empty or not. The first step towards the decision of the emptiness of the interior of $\mathcal{M}$ is the one described above, the second step being the procedure the follows from Theorem \ref{teonoiseless}.

\begin{algo}\label{algden}
	$\,$
		
	\begin{enumerate}
		\item Choose $\delta_n$ and $m_n$ in the conditions stated above;
		
		\item Compute the function $d_{\widehat{\mu},m_n}$ in the points $X_i$ using the expression (\ref{d});
		
		\item Remove from the sample the points $X_i$ where $d_{\widehat{\mu},m_n}(X_i)>\delta_n$.
		
		\item With the remaining points, use the methodology described by Theorem \ref{teonoiseless}:
		
\begin{enumerate}[label*=\arabic*]
	\item Decide $\mathring{\mathcal{M}}=\emptyset$ if and only if $peel(\hat{S}_n(r_n))=\emptyset$
\end{enumerate}

	\end{enumerate}
\end{algo}

\subsection{Simulation study}

As a way to show the results that these methods provide, we considered the manifold to be the ring with outer and inner radius of $1+\varepsilon/2$ and $1-\varepsilon/2$ respectively, that is, $\mathcal{M}=B(0,1+\varepsilon/2)\backslash B(0,1-\varepsilon/2)$, with $\varepsilon$ taking values 0,0.01,0.05 and 0.1. We drew 100 samples of size $n=250,500,1000,2500,5000$ for the case $\varepsilon=0$ (where $\mathring{\mathcal{M}}=\emptyset$) and 1000 samples of size $n=5,10,25,50,100$ for the remaining cases (where $\mathring{\mathcal{M}}\neq\emptyset$), according to the model described in the section with the following parameters:
\begin{itemize}
	\item $\mu$ as the uniform distribution on $\mathcal{M}$
	
	\item $\mu_R$ as the uniform distribution on $[-2,2]\times[-2,2]$
	
	\item $\alpha_n=n^{-y}$ with $y\in\{0.75,0.8,0.9,0.95\}$
	
	\item $m_n=n^{-0.25}$
	
	\item $\delta_n=1000n^{-0.95}$
\end{itemize}

Afterwards, we applied Algorithm \ref{algden} to each of the samples to estimate the probability of a correct interior decision by the method. The results are presented in the next tables:

\begin{table}[H]
	\centering
\begin{tabular}{|c!{\vrule width 2pt}c|c|c|c|c|}
	\hline
	y\textbackslash n & 250 & 500 & 1000 & 2500 & 5000 \\
	\ChangeRT{2pt}
	0.75 & 0.05 & 0.02 & 0.03 & 0.28 & 0.86 \\
	\hline
	0.8 & 0.05 & 0.07 & 0.05 & 0.34 & 0.94 \\
	\hline
	0.9 & 0.30 & 0.18 & 0.24 & 0.63 & 0.95 \\
	\hline
	0.95 & 0.47 & 0.34 & 0.36 & 0.75 & 0.99 \\
	\hline
\end{tabular}
\quad
\begin{tabular}{|c!{\vrule width 2pt}c|c|c|c|c|}
	\hline
	y\textbackslash n & 5 & 10 & 25 & 50 & 100 \\
	\ChangeRT{2pt}
	0.75 & 0.046 & 0.510 & 0.936 & 0.983 & 0.999 \\
	\hline
	0.8 & 0.040 & 0.475 & 0.923 & 0.961 & 0.997 \\
	\hline
	0.9 & 0.029 & 0.416 & 0.850 & 0.922 & 0.998 \\
	\hline
	0.95 & 0.034 & 0.377 & 0.851 & 0.899 & 0.997 \\
	\hline
\end{tabular}
\caption{Empirical probabilities of correct decisions for $\varepsilon=0$ (left) and $\varepsilon=0.01$ (right)}
\end{table}

\begin{table}[H]
	\centering
\begin{tabular}{|c!{\vrule width 2pt}c|c|c|c|c|}
	\hline
	y\textbackslash n & 5 & 10 & 25 & 50 & 100 \\
	\ChangeRT{2pt}
	0.75 & 0.047 & 0.583 & 0.990 & 0.999 & 1 \\
	\hline
	0.8 & 0.038 & 0.552 & 0.977 & 0.999 & 1 \\
	\hline
	0.9 & 0.040 & 0.509 & 0.966 & 0.999 & 1 \\
	\hline
	0.95 & 0.036 & 0.470 & 0.970 & 0.997 & 1 \\
	\hline
\end{tabular}
\quad
\begin{tabular}{|c!{\vrule width 2pt}c|c|c|c|c|}
	\hline
	y\textbackslash n & 5 & 10 & 25 & 50 & 100 \\
	\ChangeRT{2pt}
	0.75 & 0.049 & 0.643 & 0.997 & 1 & 1 \\
	\hline
	0.8 & 0.043 & 0.624 & 0.990 & 1 & 1 \\
	\hline
	0.9 & 0.037 & 0.581 & 0.988 & 1 & 1 \\
	\hline
	0.95 & 0.036 & 0.535 & 0.993 & 1 & 1 \\
	\hline
\end{tabular}
\caption{Empirical probabilities of correct decisions for $\varepsilon=0.05$ (left) and $\varepsilon=0.1$ (right)}
\end{table}

\end{document}